\documentclass[12pt]{article}
\usepackage{epsfig}

\newtheorem{prop}{Proposition}
\newtheorem{alg}{Algorithm}
\newtheorem{example}{Example}

\title{Numerical Evaluation of Generalized Hypergeometric Functions
for Degenerated Values of Parameters}
\author{Yasushi Tamura\\tamura@math.kobe-u.ac.jp}
\IfFileExists{graphicx.sty}{\usepackage{graphicx}}{}
\IfFileExists{epsfig.sty}{\usepackage{epsfig}}{}

\begin{document}

\maketitle

{\Large \bf Introduction} \\

In this paper,
we give an algorithm to generate connection formulas
of generalized hypergeometric functions ${}_p F _{p-1}$ 
for degenerated values of parameters.
We also show that these connection formulas give a fast method
for numerical evaluation of generalized hypergeometric functions
near $\infty$.

Several methods to evaluate generalized hypergeometric functions
are known; see, e.g., the famous text book
``Numerical Recipes''\cite{Numerical-chap4}.
As Van Der Hoeven proved, evaluating series gives a fast method
when the precision is big.
Since we need high precision values of generalized hypergeometric functions,
we will use series expansions of generalized hypergeometric
functions near $\infty$ for numerical evaluation.
We call the series expansion of ${}_p F _{p-1}(z)$ at $z = \infty$
a connection formula.
Several methods to obtain connection formulas in degenerate cases are
known among experts, but
an algorithmic method which is fast and relevant for numerical evaluation
is not known.
Van Der Hoeven gave a method to construct series solutions
around regular singular points by introducing an order among
logarithmic monomials $x^m (\log x)^n$ \cite{VanDerHoeven}.
Note that Saito, Sturmfels and Takayama found an analogous method
which is generalized to several variable case \cite[Chapter 2]{SST-ch4}.
By utilizing these methods, we will give a new method to obtain connection formulas
in degenerate cases.



Methods discussed in this paper are used in our numerical checker
for a digital formula book project \cite{tfb} and are implemented in
by Risa/Asir.
All timing data are taken by Risa/Asir on a machine 
with the following specification;
AMD Athlon MP 1800+ 1533.40-MHz, memory 2Gb PC/AT machine.

\section{Hypergeometric Function}

The function defined by the following series
is called the Gauss hypergeometric function.
\[
	F(\alpha, \beta, \gamma; z)
	= \sum_{k=0}^\infty \frac{(\alpha)_k(\beta)_k}{(\gamma)_k(1)_k} z^k
\]
\[	(a)_k = (a) (a+1) \cdots (a+k-1)	\]
The analytic continuation of this function is 
also called the Gauss hypergeometric function.
This function satisfies the following differential equation.
\[
	\{ \delta_z (\delta_z + \gamma -1) - z (\delta_z + \alpha)(\delta_z + \beta) \} F(\alpha, \beta, \gamma; z) = 0
\]
\[
	\delta_z = z \frac{\partial}{\partial z}
\]
This differential equation is called
the Gauss hypergeometric differential equation.
By expanding the products $\delta_z$ in terms of $z$ and 
$\frac{\partial}{\partial z}$, we obtain
\[
	\{ z(z-1) \frac{\partial^2}{\partial z^2}
	+ (\gamma - (\alpha+\beta+1)z ) \frac{\partial}{\partial z}
	- \alpha \beta \} F(\alpha, \beta, \gamma; z) = 0
\]
Hence, the Gauss hypergeometric function has singularities 
at $z = 0, 1, \infty$.
It is known that the hypergeometric function 
have the following integral representation (Euler integral
representation)
\[
	F(\alpha, \beta, \gamma; z)
	= \frac{\Gamma(\gamma)}{\Gamma(\alpha)\Gamma(\gamma-\alpha)}
	\int_0^1 t^{\alpha-1} (1-t)^{\gamma-\alpha-1} (1-tz)^{-\beta} dt
\]

Our goal is numerical evaluation of 
the generalized hypergeometric function
defined by the following series
\[
	_p F_{p-1}(\alpha_1, \cdots, \alpha_p,\beta_1, \cdots, \beta; z)
	= \sum_{k=0}^\infty
	\frac{(\alpha_1)_k\cdots(\alpha_p)_k}{(\beta_1)_k\cdots(\beta_{p-1})_k(1)_k}
	z^k
\]

\noindent
The famous book ``Numerical recipes'' \cite{Numerical-chap4} says that
``a fast, general routine for the complex hypergeometric function
${}_2 F_1(a,b,c;z)$, is difficult or impossible''.
One difficulty is that 
the generalized hypergeometric series 
converges only in $|z| < 1$,
then the method of evaluating series can be used only when
$|z| < 1$.
However, in case of the generalized hypergeometric functions,
there are connection formulas,
by which we can express the hypergeometric series ${}_p F_{p-1}$
in terms of a set of series which converges at $z=\infty$.
For example, in case of ${}_2 F _1$, the connection formula is as follows.
\begin{prop}[connection formula, see, e.g., \cite{GaussPainleve-chap4} ]
Assume $\alpha-\beta, \beta-\alpha, \gamma
	\not\in {\bf Z}_{\leq 0}$. Then, we have
\begin{eqnarray*}
F(\alpha,\beta,\gamma;z) &=& \frac{\Gamma(\gamma)\Gamma(\beta-\alpha)}
         {\Gamma(\gamma-\alpha)\Gamma(\beta)}
     F(\alpha, \alpha-\gamma+1, \alpha-\beta+1) (-z)^{-\alpha} \\
&+& \frac{\Gamma(\gamma)\Gamma(\alpha-\beta)}
         {\Gamma(\gamma-\beta)\Gamma(\alpha)}
     F(\beta, \beta-\gamma+1, \beta-\alpha+1) (-z)^{-\beta}
\end{eqnarray*}
\end{prop}
The condition 
 $\alpha-\beta, \beta-\alpha, \gamma
	\not\in {\bf Z}_{\leq 0}$ means
that $ \alpha-\beta, \beta-\alpha, \gamma $ are not
in the set $\{0, -1, -2, -3, \ldots \}$.
When the condition is satisfied, we say that the parameters are generic
and when the condition is not satisfied, we say that the parameters are non-generic
or degenerated.

We can use the connection formula to evaluate numerically the hypergeometric function
near $z=\infty$ in generic case.
Degenerated cases will be discussed in the rest of this paper.
We  note that the evaluating hypergeometric series can be accelerated by
the binary splitting algorithm \cite{BSA}, \cite{VanDerHoeven}.

\section{Connection Formulas of $_p F_{p-1}$}

In this section, we study connection formulas 
of $_p F_{p-1}$ between $0$ and $\infty$ and numerical evaluation
by using the formula.
In the case that parameters are generic, these formulas are
well-known. They can be obtained by using the Barns integral representation
(see, e.g., \cite[Chapter2 4.6]{GaussPainleve-chap4}).
When parameters are non-generic, there is no complete list of
connection formulas nor an algorithm to obtain connection formulas
in the literatures. 
We will give an algorithm to derive connection formulas
when parameters are non-generic.
The Algorithms 1 and 2 seem to be implicitly used among experts to study
global behaviors of generalized hypergeometric functions,
but Algorithm 3 will be new and it gives a fast routine.

Let us review the connection formula in the generic case.
It will be the starting data to generate connection formulas
in non-generic case.
\begin{prop}[connection formulas of ${}_p F_{p-1}$, 
see, e.g., \cite{GaussPainleve-chap4}]
Assume $\alpha_i-\alpha_j(i \not= j) \not\in {\bf Z}, \beta_i \not\in {\bf Z}_{\leq 0}$.
Then, we have \\
$	_p F_{p-1}(\alpha_1,\cdots,\alpha_p,\beta_1,\cdots,\beta_{p-1};z) \\
	= \sum_{i=1}^{p}
	\Pi_{j=1}^{p-1} \frac{\Gamma(\beta_j)}{\Gamma(\beta_j-\alpha_i)}
	\Pi_{j=1,i \not=j}^{p} \frac{\Gamma(\alpha_j-\alpha_i)}{\Gamma(\alpha_j)}
	{}_p F_{p-1} (\alpha_i,\alpha_i-\beta_1+1,\cdots,\alpha_i-\beta_{p-1}+1, \\
		\alpha_i-\alpha_1+1,\cdots,\alpha_i-\alpha_{i-1}+1,\alpha_i-\alpha_{i+1}+1,\cdots,\alpha_i-\alpha_p+1;
		1/z) (-z)^{-\alpha_i}	$
\end{prop}

\begin{alg}  \rm
The case of $\alpha_1 -\alpha_2 \in {\bf Z}$:

\begin{enumerate}
\item	Use contiguity relations to make $\alpha_1 = \alpha_2$.
\item	Multiply $(\alpha_2-\alpha_1)$ to the both sides
        of the connection formula in the generic case.
\item	Replace $(\alpha_2-\alpha_1)\Gamma(-\alpha_1+\alpha_2)$
by $\Gamma(-\alpha_1+\alpha_2 +1)$ and
$(\alpha_2-\alpha_1)\Gamma(\alpha_1-\alpha_2)$ by 
$-\Gamma(\alpha_1-\alpha_2 +1)$.
\item	Apply $\frac{\partial}{\partial \alpha_2}$ for the both sides
and take the limit $\alpha_2 \rightarrow \alpha_1$.
\end{enumerate}
\end{alg}


This method can be generalized to more degenerated case.

\begin{alg} \rm \label{algorithm-degenerated-p}
We assume $\alpha_1 = \alpha_2 = \cdots = \alpha_q$ $(q \leq p)$
and $^\forall i, ^\forall j, \alpha_i - \alpha_j \not\in {\bf Z}$
$(q < i, j \leq p, i \not= j)$:
\begin{enumerate}
\item	Multiply $\Pi_{i=1}^{q-1} \Pi_{j=i+1}^{q} (\alpha_j-\alpha_i)$
to the both sides of the connection formula in the generic case.
\item	Replace
$(\alpha_j-\alpha_i)\Gamma(-\alpha_i+\alpha_j)$ by 
$\Gamma(-\alpha_i+\alpha_j +1)$ and
$(\alpha_j-\alpha_i)\Gamma(-\alpha_j+\alpha_i)$ by 
$-\Gamma(-\alpha_j+\alpha_i +1)$.
\item	Apply $\frac{\partial^{\frac{q(q-1)}{2}}}{\Pi_{i=1}^{q} \partial \alpha_i^{i-1}}$ for the both sides and take the limit
$\alpha_2 \rightarrow \alpha_1, \cdots, \alpha_q \rightarrow \alpha_1$.
\end{enumerate}
\end{alg}

Note that the left hand side of the output of the algorithm is
\[	\left( \Pi_{i=1}^{q-1} i! \right) {}_p F_{p-1}(\alpha_1, \ldots, \alpha_1, 
   \alpha_{q+1},\ldots, \alpha_p, \beta_1, \ldots, \beta_{p-1};z) \]
and hence the right hand side of the output gives
a series expansion of this function around $z=\infty$.

We conjecture that
a repetition of applying the Algorithm 
\ref{algorithm-degenerated-p} and of applying contiguity relations
yields a connection formula for any degenerated case.

This method requires a complicated symbolic differentiation,
but computer algebra systems are good at it.
One can say that our method of using connection formula of ${}_p F_{p-1}$
to evaluate numerical values around $z=\infty$ is
a hybrid computation of symbolic computation and numerical computation.


\begin{example} \rm
The following is a connection formula in a degenerated case
which is obtained by the Algorithm \ref{algorithm-degenerated-p} 
\begin{eqnarray*}
	&&{}_2 F_1(\alpha_1,\alpha_1,\beta_1;z)	\\
	&=& \frac{\Gamma(\beta_1) (-z)^{-\alpha_1}}{\Gamma(\alpha_1) \Gamma(-\alpha_1 + \beta_1)}
		\lim_{\alpha_2 \rightarrow \alpha_1} \frac{\partial}{\partial \alpha_2}
			{}_2 F_1(\alpha_1, 1 + \alpha_1 - \beta_1, 1 + \alpha_1 - \alpha_2, 1/z)	\\
	&& - \frac{\Gamma(\beta_1) (-z)^{-\alpha_1}}{\Gamma(\alpha_1) \Gamma(-\alpha_1 + \beta_1)}
		\lim_{\alpha_2 \rightarrow \alpha_1} \frac{\partial}{\partial \alpha_2}
		{}_2 F_1(\alpha_2, 1 + \alpha_2 - \beta_1, 1 - \alpha_1 + \alpha_2, 1/z)	\\
	&& - \frac{\gamma   \Gamma(\beta_1) +\Gamma(\beta_1) \psi(\alpha_1)}
		{\Gamma(\alpha_1) \Gamma(-\alpha_1+\beta_1)} (-z)^{-\alpha_1}
		{}_2 F_1(\alpha_1, 1 + \alpha_1 - \beta_1, 1, 1/z)	\\
	&& - \frac{\gamma   \Gamma(\beta_1) - \Gamma(\beta_1) \log(-z)
		+ \Gamma(\beta_1) \psi(-\alpha_1 + \beta_1)}{\Gamma(\alpha_1) \Gamma(-\alpha_1 + \beta_1)}
		(-z)^{-\alpha_1}
		{}_2 F_1(\alpha_1, 1 + \alpha_1 - \beta_1, 1, 1/z).
\end{eqnarray*}
Here, $\psi(z)$ is the derivative of $\log(\Gamma(z))$
and $\gamma$ is $-\psi(1)$.
\end{example}

Our computer experiments show that 
Algorithm \ref{algorithm-degenerated-p} is not efficient.
The next algorithm is an efficient version of deriving 
connection formulas in the degenerate case.
Algorithm \ref{algorithm-degenerated-p} is used only to get $c_j^0$'s
in the following algorithm.
\begin{alg}  \label{algorithm-degenerated-p-improved}
In the case of $\alpha_1 = \alpha_2 = \cdots = \alpha_q (\mbox{where} q=p)$,
series solutions at $z=\infty$
can be written as  
$$_p F_{p-1} = (-z)^{-\alpha_1}
	\sum_{i=0}^{\infty} \{ c_0^i + c_1^i \log(-z) + \cdots + c_{q-1}^i \log^{q-1} (-z) \} z^{-i}$$
\begin{enumerate}
\item	Derive recurrence relations for $c_j^i$ with respect to $i$
  by a generalized hypergeometric differential equation
  (see, e.g., \cite{VanDerHoeven}, \cite{SST-ch4}).
  $c_j^i$ ($i>0$) are determined by $c_j^0$'s.
\item	Obtain $c_j^0$ by applying a part of Algorithm \ref{algorithm-degenerated-p}.
\item	Get other coefficients $c_j^i$
  by using the recurrence relations for $c_j^i$ with respect to $i$.
\end{enumerate}
\end{alg}


\section{Examples of Using Connection Formulas in Degenerate Cases}

Mathematica implementation for numerical evaluations of
generalized hypergeometric functions is known to be very nice.
Here we compare numerical evaluation by our formula (derived by
our Algorithm 3), evaluation by Mathematica
and evaluation by numerical integration.

\begin{example} \rm
We try to evaluate the value of ${}_2 F_1(10/3,10/3,7/2,13+13\sqrt{-1})$.
We will present timing data of evaluation
by the series expansion derived by Algorithm 3
and Mathematica.
We present timing data on Mathematica only for reader's convinience;
it is nonsense to compare our timing data by our algorithms
and those by Mathematica,
because an algorithm used by Mathematica is not known
and implementations are done on different languages.


\noindent
{\tiny
\vskip 10pt plus 1cm minus 1pt
\begin{tabular}{|l|l|l|} \hline
terms	& time	& value	\\ \hline
5	& 0.0019sec	& 0.00004646545068423618485+0.00009888637683654298440$\sqrt{-1}$	\\ \hline
10	& 0.002762sec	& 0.00004646537447334307802263261624	\\
	&	& +0.00009888640350652418659794640828$\sqrt{-1}$	\\ \hline
20	& 0.004175sec	& 0.000046465374473393490391242220236585714989 \\
	&	& +0.000098886403506421825123991664023061171848$\sqrt{-1}$	\\ \hline
40	& 0.008503sec	& 0.0000464653744733934903912421386572707301458850337603660133374 \\
	&	& +0.0000988864035064218251239916232587199904578128942387359317473$\sqrt{-1}$ \\ \hline
80	& 0.02377sec	& 0.0000464653744733934903912421386572707301458850337603824784541 \\
	&	& +0.0000988864035064218251239916232587199904578128942387442282741$\sqrt{-1}$ \\ \hline
\end{tabular}
\vskip 10pt plus 1cm minus 1pt
}

Here, terms mean the truncation degree of hypergeometric series. 
The precision is set to $\mbox{terms} + 10$.
This timing data includes time to set $c^i_j$ in Algorithm 3.
As we see in the table, the computation is done in less than 0.03 second.

Timing data by Mathematica 4.0:

\noindent
{\tiny
\vskip 10pt plus 1cm minus 1pt
\begin{tabular}{|l|l|l|} \hline
precision	& time	& value	\\ \hline
10	& 0.04 Second	& 0.0000464654 + 0.0000988864$\sqrt{-1}$ \\ \hline
25	& 0.04 Second	& 0.00004646537447339349039124214 \\
	&	& +0.00009888640350642182512399162$\sqrt{-1}$ \\ \hline
50	& 0.07 Second	& 0.000046465374473393490391242138657270730145885033760382 \\
	&	& + 0.000098886403506421825123991623258719990457812894238744$\sqrt{-1}$ \\ \hline
\end{tabular}
\vskip 10pt plus 1cm minus 1pt
}

We also try to evaluate the same value by using
Barns integral representation and adaptive rule.

Timing Data:

\noindent
{\small
\vskip 10pt plus 1cm minus 1pt
\begin{tabular}{|l|l|l|} \hline
precition	& time	& value	\\ \hline
4	& 0.3885sec	& 0.00004509188619331184102+0.00009303289997220799852$\sqrt{-1}$ \\ \hline
5	& 0.466sec	& 0.00004508711196227677087+0.00009302906708707449476$\sqrt{-1}$ \\ \hline
6	& 0.702sec	& 0.00004508788940142053260+0.00009302899442986751422$\sqrt{-1}$ \\ \hline
\end{tabular}
\vskip 10pt plus 1cm minus 1pt
}

\vskip 10pt plus 1cm minus 1pt
\begin{figure}
\centering
\includegraphics{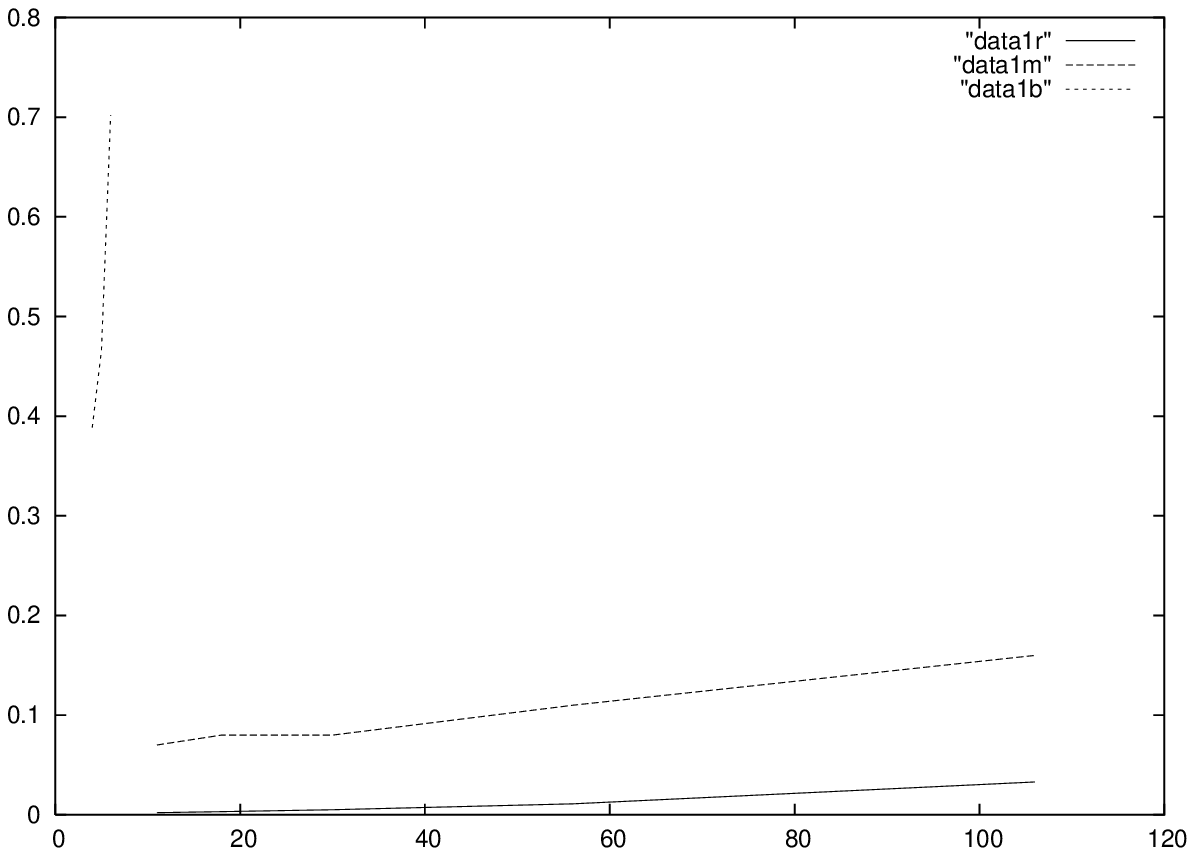}
\caption{a timing comparison of different algorithms}
\begin{tabular}{ll}
	``data1r''	&	: Algorithm 3 \\
	``data1m''	&	: Mathematica \\
	``data1b''	&	: Euler's integral representation and trapezoidal rule \\
\end{tabular}
\end{figure}
\vskip 10pt plus 1cm minus 1pt

\end{example}

\begin{example} \rm
We try to evaluate the value of
${}_2 F_1(7/2,7/2,31/5,1.3+1.8\sqrt{-1})$,
which will be more difficult than the previous example
for series expansion, since $1.3+1.8\sqrt{-1}$ is closer to the
boundary of the domain of convergence.

\noindent
{\tiny
\vskip 10pt plus 1cm minus 1pt
\begin{tabular}{|l|l|l|} \hline
terms	& time	& value	\\ \hline
5	& 0.001905sec	& -0.3879786816479458591-0.2767543538460170368$\sqrt{-1}$ \\ \hline
10	& 0.00272sec	& -0.3770255218705445491-0.2823972087891714305$\sqrt{-1}$ \\ \hline
20	& 0.004093sec	& -0.3769544095052939938707251207-0.2822863971357611098403957229$\sqrt{-1}$ \\ \hline
40	& 0.008422sec	& -0.37695442761307946514230306490910664462 \\
	&	& -0.28228642179392542114229797454872838012$\sqrt{-1}$ \\ \hline
80	& 0.02351sec	& -0.376954427613081226577499361640669979083664967552834206588 \\
	&	& -0.282286421793927502415734929810558926710399341428368162640$\sqrt{-1}$ \\ \hline
\end{tabular}
\vskip 10pt plus 1cm minus 1pt
}

The precision is set to $\mbox{terms} + 10$.
The convergence is slower than the case of 
$z = 13+13\sqrt{-1}$.
We also try to evaluate the same value by using
Euler's integral representation and trapezoidal rule.

\noindent
{\small
\vskip 10pt plus 1cm minus 1pt
\begin{tabular}{|l|l|l|} \hline
sample size	& time	& value	\\ \hline
2000	& 13.11sec	& -0.3769544276724114820-0.2822864217813210745$\sqrt{-1}$ \\ \hline
4000	& 85.41sec	& -0.3769544276222648211-0.2822864217919839264$\sqrt{-1}$ \\ \hline
8000	& 346.4sec	& -0.3769544276144992000-0.2822864217936281214$\sqrt{-1}$ \\ \hline
\end{tabular}
\vskip 10pt plus 1cm minus 1pt
}

The precision is set to $19$.
The numerical integration requires more CPU time, but 
the accuracy seems to be better than
the series expansion.
This data tell us that series expansion at $z=\infty$
should not be used around $|z|=1$.

Timing data by Mathematica 4.0:

\noindent
{\tiny
\vskip 10pt plus 1cm minus 1pt
\begin{tabular}{|l|l|l|} \hline
precision	& time	& value	\\ \hline
10	& 0.01 Second	& -0.376954 - 0.282286$\sqrt{-1}$ \\ \hline
25	& 0.01 Second	& -0.3769544276130812265774994 -0.2822864217939275024157349$\sqrt{-1}$ \\ \hline
50	& 0.01 Second	& -0.37695442761308122657749936166305176029442543627034 \\
	&	& -0.28228642179392750241573492983143989797240125238853$\sqrt{-1}$ \\ \hline
\end{tabular}
\vskip 10pt plus 1cm minus 1pt
}

\vskip 10pt plus 1cm minus 1pt
\begin{figure}
\centering
\includegraphics{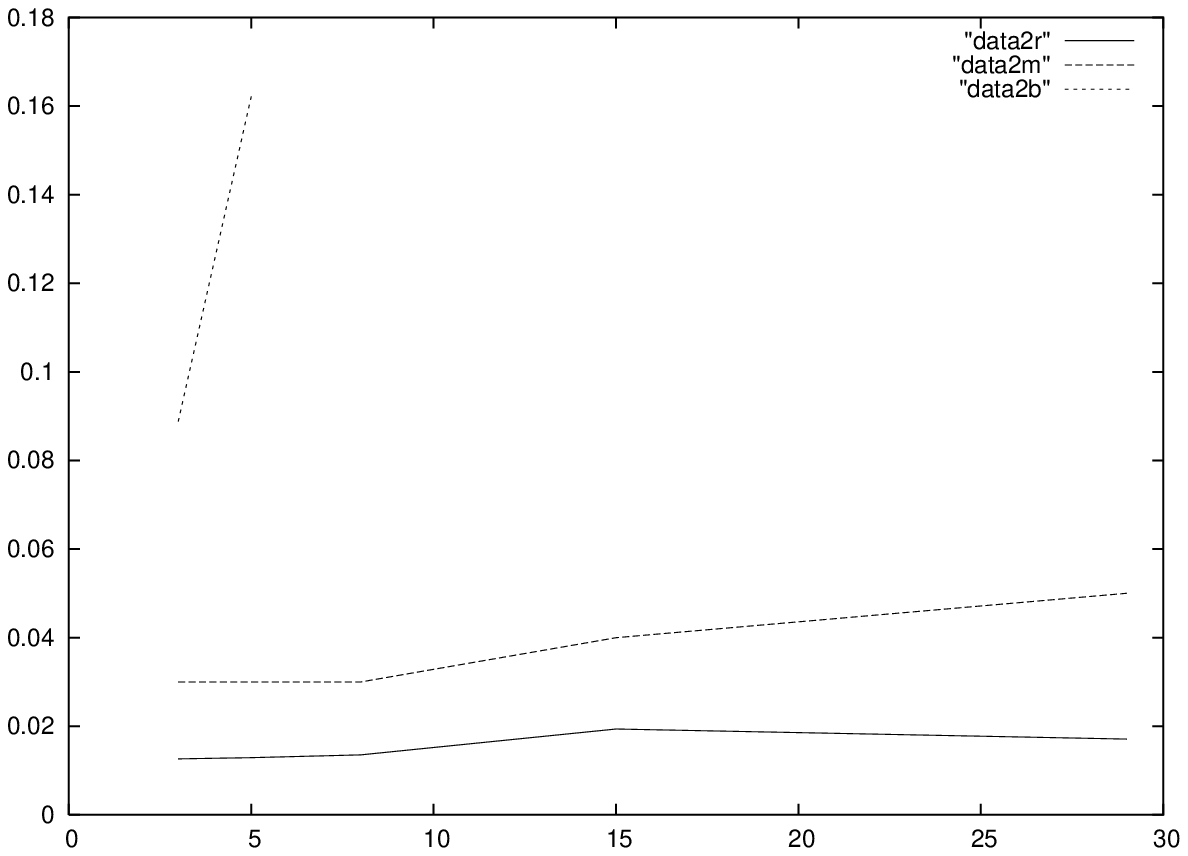}
\caption{a timing comparison of different algorithms}
\begin{tabular}{ll}
	``data2r''	&	: Algorithm 3 \\
	``data2m''	&	: Mathematica \\
	``data2b''	&	: Euler's integral representation and trapezoidal rule \\
\end{tabular}
\end{figure}
\vskip 10pt plus 1cm minus 1pt

\end{example}

The Figure below 
is a visualization
of the difference  $| {\tt r20} - {\tt r10} |$.
Here, {\tt r20} is the truncation of our connection formula
of ${}_2 F_1(7/2,7/2,31/5,x+y\sqrt{-1})$
at the degree $20$ and 
{\tt r10} is that at the degree $10$.

\vskip 10pt plus 1cm minus 1pt
\includegraphics{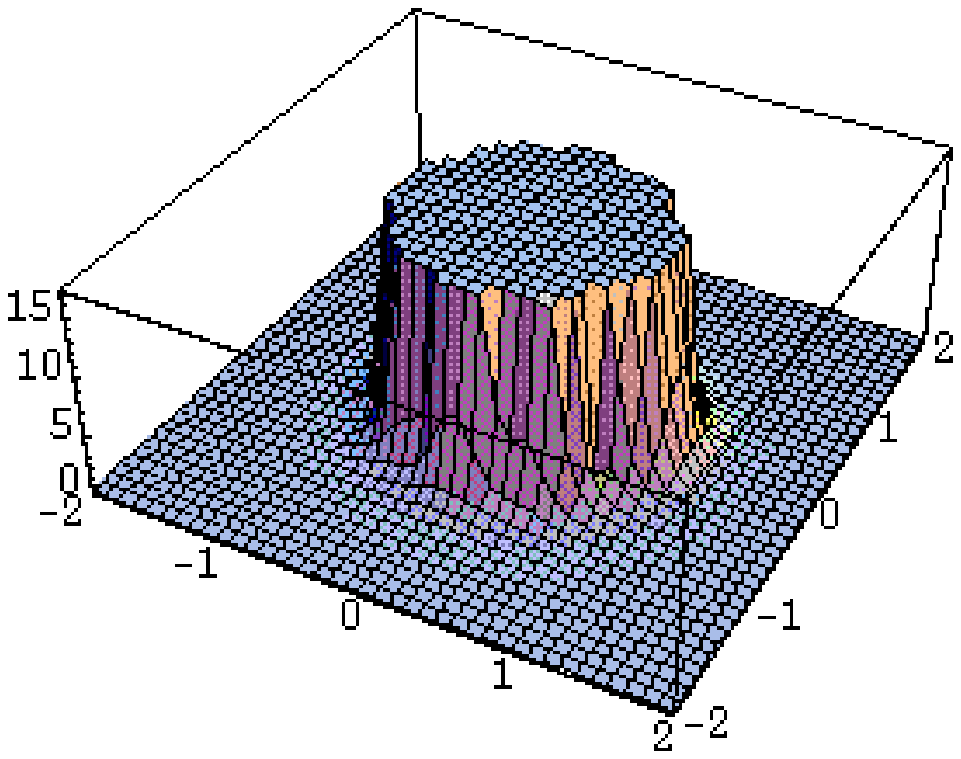}
\vskip 10pt plus 1cm minus 1pt

\noindent
The difference is close to zero in $|z| >> 1$,
but  it is larger in the neighborhood of $|z|=1$.
Developing a method for high-precision evaluation around $|z|=1$ will
be a future problem.

\begin{example} \rm
We will evaluate values of $_3 F_2(7/2,7/2,7/2,31/5,36/7;z)$.
We compare the adaptive integration method of the Barns integral
and numerical evaluation of our degenerated connection formulas
derived by Algorithm \ref{algorithm-degenerated-p-improved}.

Numerical integration(precision is set to $19$):

\noindent
{\small
\vskip 10pt plus 1cm minus 1pt
\begin{tabular}{|l|l|l|l|l|l|l|}	\hline
	& time	& value	\\ \hline
$z=130+130\sqrt{-1}$	& 0.2219sec	& 0.000001881452796232078012	\\
	&	& 0.000006893851655427774880$\sqrt{-1}$	\\ \hline
$z=13+13\sqrt{-1}$	& 0.1605sec	& 0.007312359138710618527	\\
	&	& -0.006400306129230969169$\sqrt{-1}$	\\ \hline
$z=1.3+1.3\sqrt{-1}$	& 0.112sec	& -1.097506492885595820	\\
	&	& +0.6369234717153367928$\sqrt{-1}$	\\ \hline
\end{tabular}
\vskip 10pt plus 1cm minus 1pt
}

Series expansion(terms are set to $20$, and precision is set to $19$):

\noindent
{
\vskip 10pt plus 1cm minus 1pt
\begin{tabular}{|l|l|l|l|l|l|l|}	\hline
	& time	& value	\\ \hline
$z=130+130\sqrt{-1}$	& 0.01571sec	& 0.00001345106300346753915	\\
	&	& 0.000006796099418228839164$\sqrt{-1}$	\\ \hline
$z=13+13\sqrt{-1}$	& 0.01404sec	& 0.007350815068974895610	\\
	&	& -0.006282360701607166085$\sqrt{-1}$	\\ \hline
$z=1.3+1.3\sqrt{-1}$	& 0.004741sec	& -1.097992622097576377	\\
	&	& +0.6364759787999937697$\sqrt{-1}$	\\ \hline
\end{tabular}
\vskip 10pt plus 1cm minus 1pt
}

Timing data by Mathematica 4.0(precision is set to $10$):

\noindent
\vskip 10pt plus 1cm minus 1pt
\begin{tabular}{|l|l|l|l|l|l|l|}	\hline
	& time	& value	\\ \hline
$z=130+130\sqrt{-1}$	& 0.1 Second	& 0.0000134511	\\
	&	& +6.79615 $10^{-6} \sqrt{-1}$ 	\\ \hline
$z=13+13\sqrt{-1}$	& 0.07 Second	& 0.00735089	\\
	&	& - 0.00628229$\sqrt{-1}$	\\ \hline
$z=1.3+1.3\sqrt{-1}$	& 0.07 Second	& -1.09725	\\
	&	& +0.636973$\sqrt{-1}$	\\ \hline
\end{tabular}
\vskip 10pt plus 1cm minus 1pt
\end{example}

We conclude that our formulas derived by algorithm 3 give a fast method
to evaluate generalized hypergeometric functions around $\infty$.

\end{document}